\documentclass[11pt]{amsart}
\usepackage{xcolor}

\def\R{{\mathbb R}}
\def\Z{{\mathbb Z}}
\def\C{{\mathbb C}}
\def\Pr{{\mathbb P}}
\def\Exp{{\mathbb E}}

\def\Rad{{\mathcal R}}

\def\det{\rm det}

\def\N{{\mathbb N}}

\def\e{\varepsilon}

\numberwithin{equation}{section}

\begin{document}
\title[]{Inequalities for sections and projections of convex bodies}
\author{Apostolos Giannopoulos}
\address[Apostolos Giannopoulos]{Department of Mathematics, National and Kapodistrian University of Athens,
Panepistimiopolis 157-84, Athens, Greece}
\email{apgiannop@@math.uoa.gr}

\author{Alexander  Koldobsky}
\address[Alexander  Koldobsky]{Department of Mathematics, University of Missouri, Columbia, MO 65211, USA}
\email{koldobskiya@@missouri.edu}

\author{Artem Zvavitch}
\address[Artem Zvavitch]{Department of Mathematical Sciences\\ Kent State University\\ Kent, OH USA}
\email{zvavitch@@math.kent.edu}

\thanks{The first named author is supported by the Hellenic Foundation for
Research and Innovation (H.F.R.I.) under the ``First Call for H.F.R.I. Research Projects to support Faculty members and Researchers and
the procurement of high-cost research equipment grant" (Project Number: 1849).
The second named author was supported in part by the U.S. National Science Foundation Grant DMS-2054068.
The third named author was supported in part by  the U.S. National Science Foundation Grant DMS-2000304 and  United States - Israel Binational Science Foundation (BSF). Both the second and the third named authors were  supported in part by U.S.
National Science Foundation under Grant No. DMS-1929284 while
 in residence at the Institute for
Computational and Experimental Research in Mathematics in Providence, RI, during the Harmonic Analysis
and Convexity semester program.}
\keywords{Convex bodies; Sections; Radon transform; Intersection body}
\subjclass[2010]{52A20, 53A15, 52B10.}

\begin{abstract} This article belongs to the area of geometric tomography, which is the study of geometric properties of solids based
on data about their sections and projections. We describe a new direction in geometric tomography where different volumetric results are
considered in a more general setting, with volume replaced by an arbitrary measure. Surprisingly, such a general approach works for
a number of volumetric results. In particular, we discuss the Busemann-Petty problem on sections of convex bodies
for arbitrary measures and the slicing problem for arbitrary measures.  We present generalizations of these questions to the case of functions.
A number of generalizations of questions related to projections, such as the problem of Shephard, are also discussed as well as some questions
in discrete  tomography.
\end{abstract}

\maketitle

\tableofcontents

\section{Introduction}
The Busemann-Petty problem asks whether origin-symmetric convex bodies in $\R^n$ with uniformly smaller $(n-1)$-dimensional volume of central hyperplane sections necessarily have smaller $n$-dimensional volume. The slicing problem of Bourgain asks whether every symmetric convex body of volume one in $\R^n$ has a central hyperplane section whose $(n-1)$-dimensional volume is greater than an absolute constant. We look at these and other results and problems of convex geometry from a more general point of view, replacing volume by an arbitrary measure. Though common sense suggests
that the setting of arbitrary measures is too general to produce significant results, we present several situations where such generalizations are very much possible. In particular, it was shown in \cite{Zvavitch-2005} that the solution of the Busemann-Petty problem (affirmative if $n\le 4$ and negative if $n\ge 5$) is exactly the same for an arbitrary measure with positive density in place of volume. A version of the slicing problem for arbitrary measures was proved in \cite{Koldobsky-2015b}, namely for any probability density $f$ on an origin-symmetric convex body $K$ of volume one, there exists a hyperplane $H$ in $\R^n$ so that the integral of $f$ over $K\cap H$ is greater than $\frac{1}{2\sqrt{n}}.$

The paper is organized as follows. In Section \ref{ND} we briefly introduce the most essential basic notation and facts required for our exposition. In Section \ref{section:LW-M}, first we present volume estimates from orthogonal projections and sections and continue with generalizations and
variants of these inequalities. In Section \ref{sections} we discuss the Busemann-Petty problem and its generalizations with the
emphasis on the general setting with measures in place of volume. Section \ref{projections} covers results related to projections of convex bodies.
These include the Shephard's problem which is the projection analogue of the Busemann-Petty problem, Milman's problem which can be considered as a mixed Busemann-Petty-Shephard's problem, and slicing-type inequalities for the surface area of projections. Section \ref{surface-area} deals with comparison and slicing inequalities for the surface area of convex bodies. In Section \ref{vol-diff} we discuss volume difference inequalities which allow to estimate the error in tomographic calculations. Finally, in Section \ref{discrete} we present what is known about discrete analogues of the slicing problem.

\vspace{.1in}

\noindent {\bf Acknowledgments.} We are grateful to Dylan Langharst and Michael Roysdon for many corrections, valuable discussions and suggestions.

\section{Notation and definitions}\label{ND}

In this section we will introduce a few basic notations and definitions needed for this article, we refer the reader to \cite{AGA-book, AGA-bookII, BGVV-book,  Burago-Zalgaller-book, Gardner-book, Gardner-bookadd, Koldobsky-book, Koldobsky-Yaskin-book, Schneider-book} for a wealth of additional information on objects and tools from Convex Geometry, Geometric Tomography and Fourier Analysis  used in this survey.

We work in ${\mathbb R}^n$, which is equipped with the standard inner product $\langle\cdot ,\cdot\rangle $. We denote by
$B_2^n$ and $S^{n-1}$ the Euclidean unit ball and sphere respectively. We write $|\cdot |$ for volume in the
appropriate dimension, $\omega_n$ for the volume of $B_2^n$ and $\sigma $ for the rotationally invariant probability
measure on $S^{n-1}$. The Grassmann manifold $G_{n,k}$ of all $k$-dimensional subspaces of ${\mathbb R}^n$ is
equipped with the Haar probability measure $\nu_{n,k}$. For every $1\leq k\leq n-1$ and $H\in G_{n,k}$ we denote by $P_H$
the orthogonal projection from $\mathbb R^{n}$ onto $H$. The letters $c,c^{\prime }, c_1, c_2$ etc. denote absolute positive
constants which may change from line to line. Whenever we write $a\approx b$, we mean that there exist absolute constants
$c_1,c_2>0$ such that $c_1a\leq b\leq c_2a$.

A convex body in ${\mathbb R}^n$ is a compact convex subset $K$ of ${\mathbb R}^n$ with non-empty interior. We say that $K$ is
origin-symmetric if $-K= K$, and that $K$ is centered if its barycenter $\frac{1}{|K|}\int_Kx\,dx $ is at the origin.
The support function of a convex body $K$ is defined by $h_K(y)=\max \{\langle x,y\rangle:x\in K\}$, and the mean width of $K$ is
\begin{equation*}w(K)=\int_{S^{n-1}}h_K(\xi )\,d\sigma (\xi ). \end{equation*}
A closed bounded set $K$ in $\R^n$ is called a star body if every straight line passing through the origin crosses the boundary of $K$
at exactly two points different from the origin, the origin is an interior point of $K,$ and the Minkowski functional of $K$ defined by
$$\|x\|_K = \min\{a\ge 0:\ x\in aK\}$$
is a continuous function on $\R^n.$  We use the polar formula for the volume $|K|$ of a star body $K:$
\begin{equation}\label{polar-vol}
|K|=\frac 1n \int_{S^{n-1}} \|\xi\|_K^{-n} d\xi.
\end{equation}
If $f$ is an integrable function on $K$, then
\begin{equation}\label{polar-meas}
\int_K f = \int_{S^{n-1}} \left(\int_0^{\|\xi\|_K^{-1}} r^{n-1}f(r\xi) dr \right) d\xi.
\end{equation}
For $1\le k \le n-1,$  the $(n-k)$-dimensional spherical Radon transform $\Rad_{n-k}:C(S^{n-1})\to C(G_{n,n-k})$ is a linear operator defined by
$$\Rad_{n-k}g (H)=\int_{S^{n-1}\cap H} g(x)\ dx\quad \mbox{for all}\;\;  H\in G_{n,n-k}$$
for every function $g\in C(S^{n-1}).$

\smallbreak

For every $H\in G_{n,n-k},$ the $(n-k)$-dimensional volume of the section of a star body $K$ by $H$ can be written as
\begin{equation}\label{vol-sect}
|K\cap H| = \frac1{n-k} \Rad_{n-k}(\|\cdot\|_K^{-n+k})(H).
\end{equation}
More generally, for an integrable function $f$ and any $H\in G_{n,n-k}$,
\begin{equation}\label{meas-sect}
\int_{K\cap H} f = \Rad_{n-k}\left(\int_0^{\|\cdot\|_K^{-1}} r^{n-k-1}f(r\ \cdot)\ dr \right)(H).
\end{equation}
The class of intersection bodies ${\mathcal{I}}_n$ was introduced by Lutwak \cite{Lutwak-1988}. We consider a generalization of this concept due to Zhang \cite{Zhang-1996}. We say that an origin symmetric star body $D$ in $\R^n$ is a generalized $k$-intersection body, and write $D\in {\mathcal{BP}}_k^n,$  if there exists a finite Borel non-negative measure $\nu_D$ on $G_{n,n-k}$ so that for every $g\in C(S^{n-1})$
\begin{equation}\label{genint}
\int_{S^{n-1}} \|x\|_D^{-k} g(x)\ dx=\int_{G_{n,n-k}} R_{n-k}g(H)\ d\nu_D(H).
\end{equation}
When $k=1$ we get the original Lutwak's class of intersection bodies ${\mathcal{BP}}_1^n={\mathcal{I}}_n$.

Let ${\mathcal A}$ be a class of star bodies in ${\mathbb R}^n$ which is invariant with respect to invertible linear transformations. We denote by
$$d_{BM}(K,{\mathcal A})=\inf \{a>0: \exists D\in {\mathcal A} \quad {\rm such\ that} \quad K \subset D \subset aK\}$$
the Banach-Mazur distance from $K$ to ${\mathcal A}$. We also define the smaller volume ratio distance
$$d_{\rm {vr}}(K,{\mathcal A}) = \inf \left\{ \left( |K|/|D|\right)^{1/n}:\ D\subset K,\ D\in {\mathcal A} \right\},$$
and outer volume ratio distance
$$d_{\rm {ovr}}(K,{\mathcal A}) = \inf \left\{ \left( |D|/|K|\right)^{1/n}:\ K\subset D,\ D\in {\mathcal A} \right\}$$
from $K$ to ${\mathcal A}$.

Minkowski's fundamental theorem states that if $K_1,\ldots ,K_m$ are non-empty, compact convex
subsets of ${\mathbb R}^n$, then the volume of $t_1K_1+\cdots +t_mK_m$ is a homogeneous polynomial of degree $n$ in
$t_i>0$; that is,
\begin{equation*}|t_1K_1+\cdots +t_mK_m|=\sum_{1\leq i_1,\ldots ,i_n\leq m}
V(K_{i_1},\ldots ,K_{i_n})t_{i_1}\cdots t_{i_n},\end{equation*}
where the coefficients $V(K_{i_1},\ldots ,K_{i_n})$ are chosen to be invariant under permutations of their arguments.
The coefficient $V(K_1,\ldots ,K_n)$ is the mixed volume of $K_1,\ldots ,K_n$; we refer to \cite{Schneider-book}
for a detailed exposition of the definition and main properties of mixed volumes. In particular, if $K$ and $D$ are two convex bodies in ${\mathbb R}^n$,
then the function $|K+tD|$ is a polynomial in $t\in [0,\infty )$:
\begin{equation*}
|K+tD|=\sum_{j=0}^n \binom{n}{j} V_{n-j}(K,D)\;t^j,
\end{equation*}
where $V_{n-j}(K,D)= V((K,n-j),(D,j))$ is the $j$-th mixed volume of $K$ and $D$ (we use  the notation $(D,j)$ for $D,\ldots ,D$ $j$-times).
If $D=B_2^n$ then we set $W_j(K):=V_{n-j}(K,B_2^n)=V((K, n-j), (B_2^n, j))$; this is the $j$-th quermassintegral of $K$.
The mixed volume $V_{n-1}(K,D)$ can be expressed as
\begin{equation}\label{eq:not-1}V_{n-1}(K,D)={\frac{1}{n}}\int_{S^{n-1}}h_D(\xi )d\sigma_K(\xi ),\end{equation}
where $\sigma_K$ is the surface area measure of $K$; this is the Borel measure
on $S^{n-1}$ defined by
\begin{equation*}\sigma_K(A)=\lambda (\{x\in {\rm bd}(K):\;{\rm the}\;{\rm outer}
\;{\rm normal}\;{\rm to}\;K\;{\rm at}\;x\;{\rm belongs}\;
{\rm to}\;A\}),\end{equation*}
where $\lambda $ is the Hausdorff measure on ${\rm bd}(K)$.  In particular, if $\sigma_K$ is absolutely continuous with respect
to $\lambda$ then the density of $\sigma_K$ is called  the curvature function and is usually denoted as $f_K$. The surface area $S(K):=\sigma_K(S^{n-1})$ of $K$ satisfies
\begin{equation*}S(K)=nW_1(K).\end{equation*}
Volume and mixed volumes in general satisfy a number of very useful inequalities. The first one is the Brunn-Minkowski inequality
$|K+L|^{1/n}\ge |K|^{1/n}+|L|^{1/n},$ whenever  $K, L$ and $K+L$ are measurable and nonempty.

Direct consequences of the Brunn-Minkowski inequality  are   Minkowski's first inequality
\begin{equation}\label{M1}
V_{n-1}(K, L)\ge |K|^{(n-1)/n}|L|^{1/n},
\end{equation}
and  Minkowski's second inequality
\begin{equation}\label{M2}
V(K, L)^2\ge |K|V((L,2),(K, n-2)),
\end{equation}
for two convex, compact subsets $K$ and  $L$ of $\R^n$.

A zonoid is the limit of Minkowski sums of line segments in the Hausdorff metric. Equivalently, an
origin-symmetric convex body $Z$ is a zonoid if and only if its polar body $Z^{\circ}$ is the unit ball of an $n$-dimensional subspace of an $L_1$-space;
i.e. if there exists a positive measure $\mu $ (the supporting measure of $Z$) on $S^{n-1}$ such that
\begin{equation*}h_Z(x)=\| x\|_{Z^{\circ }}=\frac{1}{2}\int_{S^{n-1}}|\langle x,y\rangle |d\mu (y).\end{equation*}
The class of origin-symmetric zonoids with non-empty interior coincides with the class of projection bodies.
Recall that the projection body $\Pi K$ of a convex body $K$ is the
symmetric convex body whose support function is defined by
\begin{equation*}h_{\Pi K} (\xi )=|P_{\xi^{\perp } }(K)|=\frac{1}{2}\int_{S^{n-1}}|\langle \xi ,y\rangle |d\sigma_K(y), \qquad \xi\in S^{n-1},\end{equation*}
where $\xi^\perp=\{x\in \R^n:\ \langle x,\xi\rangle =0\}$ is the central hyperplane perpendicular to $\xi$ and $P_{\xi^{\perp }}(K)$ denotes the orthogonal projection
of $K\subset {\mathbb R}^n$ onto  $\xi^{\perp }$.

\section{Volume estimates from orthogonal projections and sections}\label{section:LW-M}

Estimating the volume of a convex body from the volumes of its orthogonal projections or sections is a classical
question in convex geometry. A well-known such estimate is the famous Loomis-Whitney inequality, which asserts that for any convex body (actually, any compact set) $K$
\begin{equation}\label{LW}
|K|^{n-1} \leq \prod\limits_{i=1}^n |P_{e_i^\perp}(K)|,
\end{equation}
where $\{e_1, \ldots, e_n\}$ is an orthonormal basis of $\R^n$  (see \cite{Loomis-Whitney-1949}).
Equality holds in \eqref{LW} if and only if $K$ is an orthogonal
parallelepiped such that $\pm e_i$ are the normal vectors of its facets. A dual inequality, in which the volume of $K$ is estimated
by the volumes of the sections $K\cap e_i^{\perp }$, was proved by Meyer in \cite{Meyer-1988}: for every convex body $K$ in ${\mathbb R}^n$ one has
\begin{equation}\label{eq:LW-2}
|K|^{n-1}\geq \frac{n!}{n^n}\prod_{i=1}^n|K\cap e_i^{\perp }|,
\end{equation}
with equality if and only if $K={\rm conv}\{\pm \lambda_1e_1,\ldots ,\pm \lambda_ne_n\}$
for some $\lambda_i>0$.

Both inequalities admit various generalizations. In order to state one of them, let $s>0$ and say that the subspaces
$F_1,\ldots ,F_r$ form an $s$-uniform cover of ${\mathbb R}^n$ with weights $c_1,\ldots ,c_r>0$ if
\begin{equation}\label{eq:uniform-weights}s\,I_n=\sum_{i=1}^rc_iP_{F_i},\end{equation}
where $I_n$ is the identity operator on ${\mathbb R}^n$.
Then, as an application of the multidimensional geometric Brascamp-Lieb inequality, one may show that,
for every compact subset $K$ of ${\mathbb R}^n,$ we have
\begin{equation}\label{eq:s-cover-1}
|K|^s\leq\prod_{i=1}^r|P_{F_i}(K)|^{c_i}.
\end{equation}
On the other hand, using Barthe's geometric reverse Brascamp-Lieb inequality,
it was proved in \cite{Liakopoulos-2018} that if $K$ is a convex body in
${\mathbb R}^n$ with $0\in {\rm int}(K)$ and $F_1,\ldots ,F_r$ are subspaces as above, then
\begin{equation}\label{eq:s-cover-2}
|K|^s\geq\frac{1}{(n!)^s}\prod_{i=1}^r\big (d_i!\,|K\cap F_i|\big)^{c_i},
\end{equation}
where $d_i={\rm dim}(F_i)$. A special case of these inequalities occurs when $u_1,\ldots ,u_m$ are unit vectors in ${\mathbb R}^n$ and
$c_1,\ldots ,c_m$ are positive real numbers that satisfy John's condition
$I_n=\sum_{i=j}^mc_ju_j\otimes u_j$. Then, if $K$ is a convex body  in
${\mathbb R}^n$ with $0\in {\rm int}(K),$ we have that
\begin{equation}\label{eq:LW-3}
\frac{n!}{n^n}\prod_{j=1}^m|K\cap u_j^{\perp }|^{c_j}\leq
|K|^{n-1}\leq\prod_{j=1}^m|P_{u_j^{\perp }}(K)|^{c_j}.
\end{equation}
To see this, observe that if $P_j=P_{u_j^{\perp }}$
then $u_j\otimes u_j=I_n-P_j$; hence John's condition may be written as $I_n=\sum_{j=1}^mc_j(I_n-P_j)$, which implies that
\begin{equation}\label{eq:LW-ball-1}(n-1)I_n=\sum_{j=1}^mc_jP_j,\end{equation}
because $\sum_{j=1}^mc_j=n$. The assumption that the interior of $K$ contains the origin is needed only for the left hand side inequality.
The right hand side inequality in \eqref{eq:LW-3} was proved by Ball in \cite{Ball-1991}, while
the left hand side inequality was obtained by Li and Huang in \cite{Li-Huang-2016}.

Another extension of the Loomis-Whitney inequality, which can be put in the same framework, had been established
in \cite{Bollobas-Thomason-1995}. For every non-empty $\tau\subset [n],$ where $[n]=\{1, \dots, n\},$ we set $F_{\tau }={\rm span}\{e_j:j\in\tau \}$ and $E_{\tau }=F_{\tau }^{\perp }$.
Given an integer $s\geq 1,$ we say that the (not necessarily distinct) sets $\sigma_1,\ldots ,\sigma_r\subseteq [n]$
form an $s$-uniform cover of $[n]$ if every $j\in [n]$ belongs to exactly $s$ of the sets $\sigma_i$.
The uniform cover inequality of Bollob\'{a}s and Thomason states that, for every compact subset $K$ of ${\mathbb R}^n$ which is the closure
of its interior, we have
\begin{equation}\label{eq:LW-4}
|K|^s\leq\prod_{i=1}^r|P_{F_{\sigma_i}}(K)|.
\end{equation}
Note that if $(\sigma_1,\ldots ,\sigma_r)$ is an $s$-uniform cover
of $[n]$, then setting $F_i=F_{\sigma_i}={\rm span}(\{e_j:j\in\sigma_i\})$, $i\in [r]$, we have $s\,I_n=\sum_{i=1}^rP_{F_i}$.
Thus, \eqref{eq:LW-4} is an immediate consequence of \eqref{eq:s-cover-1}. Also, \eqref{eq:s-cover-2} implies
that if $K$ is a convex body  in ${\mathbb R}^n$ with $0\in {\rm int}(K),$ then
\begin{equation}\label{eq:dual-BT-L}
|K|^s\geq\frac{1}{(n!)^s}\prod_{i=1}^r|\sigma_i|!\,|K\cap F_i|.
\end{equation}
To recover Meyer's inequality from \eqref{eq:s-cover-2}, we use the particular case $F_i=e_i^{\perp }$, $i\in [n]$, so that
we have $(n-1)\,I_n=\sum_{i=1}^nP_{e_i^{\perp }}$. Applying \eqref{eq:dual-BT-L}
with $s=n-1$ and $|\sigma_i|={\rm dim}(F_i)=n-1$ we get
\begin{equation*}
|K|^{n-1}\geq\frac{n!}{n^n}\prod_{i=1}^n|K\cap e_i^{\perp }|
\end{equation*}
for any convex body $K$ in ${\mathbb R}^n$ with $0\in {\rm int}(K)$.
It is also not hard to see that the left hand side inequality in \eqref{eq:LW-3} is a consequence of \eqref{eq:s-cover-2}.

Local Loomis-Whitney type inequalities were studied in many  works, including \cite{Giannopoulos-Hartzoulaki-Paouris-2002,  Fradelizi-Giannopoulos-Meyer-2003}, where it was proved that for any convex body $K$ in $\R^n$ and a pair of orthogonal vectors
$u,v \in S^{n-1}$, one has
$$
|K||P_{[u,v]^\perp}(K)| \le \frac{2(n-1)}{n}  |P_{u^\perp} (K)||P_{v^\perp} (K)|.
$$
We refer to \cite{Fradelizi-Madiman-Meyer-Zvavitch2022} for a simple proof of this inequality and a number of equivalent restatements.

Many restricted variants of the Loomis-Whitney inequality and of the uniform cover inequality, estimating the volume of a convex
body from the volumes of a smaller set of sections or projections, were
obtained in \cite{BGL}. See also \cite{Soprunov-Zvavitch-2016, Alonso-Artstein-Gonzales-Jimenez-Villa-2019, Fradelizi-Madiman-Meyer-Zvavitch2022}
for some sharp results in this direction.

\section{Comparison and slicing inequalities for functions}\label{sections}

\subsection{The comparison problem for functions}

In 1956, Busemann and Petty \cite{Busemann-Petty-1956} posed the
problem if, for any origin-symmetric convex bodies $K,L$ in $\R^n$, the inequalities
\begin{equation}\label{bp-condition}
|K\cap\xi^\perp | \le |L\cap\xi^\perp|\quad \mbox{for all}\;\xi\in S^{n-1}
\end{equation}
imply $|K| \le |L|$.
The problem was solved at the end of the 1990's in a sequence of papers \cite{Larman-Rogers-1975, Ball-1988, Giannopoulos-1990, Bourgain-1991b, Lutwak-1988, Papadimitrakis-1992, Gardner-1994, Gardner-1994b, Zhang-1994, Koldobsky-1998, Koldobsky-1998b, Zhang-1999, Gardner-Koldobsky-Schlumprecht-1999}. The answer is affirmative if $n\le 4$, and it is negative if $n\ge 5.$
We refer the reader to \cite[p.~3]{Koldobsky-book} or \cite[p.~343]{Gardner-book} for the history of the solution.

The isomorphic Busemann-Petty problem, posed in \cite{VMilman-Pajor-1989}, asks whether the inequalities (\ref{bp-condition}) imply $\left|K\right| \le C\left|L\right|,$
where $C$ is an absolute constant. This question is still open and is equivalent to the slicing problem of Bourgain; see below. A recent result
of Klartag and Lehec \cite{Klartag-Lehec-2022} shows that the constant can be estimated by a polylogarithmic function of the dimension.

An extension of the Busemann-Petty problem to arbitrary measures in place of volume was considered in \cite{Zvavitch-2005}.
Let $K,L$ be origin-symmetric convex bodies in $\R^n,$ and let $f$ be a locally integrable non-negative function on $\R^n.$
Suppose that for every $\xi\in S^{n-1}$
\begin{equation}\label{bp}
\int_{K\cap \xi^\bot} f(x) dx \le \int_{L\cap \xi^\bot} f(x) dx,
\end{equation}
where integration is with respect to Lebesgue measure on $\xi^\bot.$
Does it necessarily follow that
$$
\int_K f(x) dx \le \int_L f(x) dx?
$$
It was proved in \cite{Zvavitch-2005} that, for any strictly positive function $f$, the solution is the same as in the case of volume (where $f\equiv 1$):
affirmative if $n\le 4$ and negative if $n\ge 5$.

In view of this result, it is natural to ask the isomorphic question again. Do inequalities (\ref{bp}) imply that
$$
\int_K f(x) dx \le s_n \int_L f(x) dx,
$$
where the constant $s_n$ does not depend on $f,K,L?$ It was proved in \cite{Koldobsky-Zvavitch-2015} that the answer is affirmative, namely $s_n\le \sqrt{n}.$

The argument in \cite{Koldobsky-Zvavitch-2015} is based on a more general estimate. It was proved in \cite{Koldobsky-Zvavitch-2015} that the inequalities (\ref{bp}) imply
\begin{equation}\label{bm}
\int_K f(x) dx \le d_{BM}(K,{\mathcal{I}}_n) \int_L f(x) dx.
\end{equation}
By John's theorem \cite{John-1948} and the fact that the class ${\mathcal{I}}_n$ contains ellipsoids, we have $d_{BM}(K,{\mathcal{I}}_n)\le \sqrt{n}$,
which proves the $\sqrt{n}$ estimate in the isomorphic Busemann-Petty problem for functions.
It is not known whether the $\sqrt{n}$ estimate is optimal. Another open question is whether the Banach-Mazur distance in (\ref{bm})
can be replaced by the smaller outer volume ratio distance $d_{\rm {ovr}}(K,{\mathcal{I}}_n)$.

A slightly different estimate with the outer volume ratio distance instead of the Banach-Mazur distance was proved in \cite{Koldobsky-Paouris-Zvavitch-2022}. Namely,
if $K,L$ are star bodies in $\R^n,$ and $f,g$ are non-negative locally integrable functions on $\R^n$ with $\|g\|_\infty=g(0)=1$, then the inequalities
$$\int_{K\cap \xi^\bot} f(x) dx \le \int_{L\cap \xi^\bot} g(x) dx\quad \mbox{for all}\;\; \xi\in S^{n-1}$$
imply
\begin{equation}\label{bp1}\int_K f(x) dx\le d_{\rm ovr}(K,{\mathcal{I}}_n) \frac n{n-1} |K|^{\frac 1n}\left(\int_L g(x) dx \right)^{\frac{n-1}n}.\end{equation}

\subsection{The slicing problem for functions}

The slicing problem of Bourgain \cite{Bourgain-1986, Bourgain-1987}
asks whether there exists a constant $C$ so that, for any $n\in \N$ and any
origin-symmetric convex body $K$ in $\R^n,$
\begin{equation} \label{hyper}
|K|^{\frac {n-1}n} \le C \max_{\xi \in S^{n-1}} |K\cap \xi^\bot|.
\end{equation}
In other words, is it true that every origin-symmetric convex
body $K$ in $\R^n$ of volume $1$ has a hyperplane section whose $(n-1)$-dimensional
volume is greater than an absolute constant?

The problem remains open. Bourgain \cite{Bourgain-1991} proved that $C\le O(n^{1/4})$ up to a logarithmic factor which was removed by Klartag \cite{Klartag-2006}. Chen \cite{Chen-2021} proved that $C\le O(n^{\epsilon})$ for every $\epsilon>0,$ and Klartag and Lehec \cite{Klartag-Lehec-2022} established a polylogarithmic bound $C\le O(\log^4n)$. The method of \cite{Klartag-Lehec-2022} was slightly refined in \cite{JLV} where it was shown that
$C\le O(\log^{2.2226}n).$ The answer is known to be affirmative for some special classes of convex bodies.
For unconditional convex bodies this was observed by Bourgain; see also \cite{VMilman-Pajor-1989, Junge-1995,
Bobkov-Nazarov-2003}), for unit balls of subspaces of $L_p$ it was proved in \cite{Ball-1991b, Junge-1994, EMilman-2006}, for intersection bodies in
\cite[Th.9.4.11]{Gardner-book}, for zonoids, duals of bodies with bounded volume ratio in
\cite{VMilman-Pajor-1989}, for the Schatten classes in \cite{Konig-Meyer-Pajor-1998}, and for $k$-intersection bodies in \cite{Koldobsky-Pajor-Yaskin-2008, Koldobsky-2016}.
Other partial results on the problem include \cite{Ball-1988b, Bourgain-Klartag-Milman-2004, Dafnis-Paouris-2010, Dar-1995, Giannopoulos-Paouris-Vritsiou-2012, Klartag-2005,
Klartag-Kozma-2009, Paouris-2000, Eldan-Klartag-2011, Ball-Nguyen-2012}; see the book \cite{BGVV-book} and the surveys \cite{Klartag-Milman, Gardner-bookadd} for details.

A generalization of the slicing problem to arbitrary functions was suggested in \cite{Koldobsky-2012}.
Does there exist a constant $T_n$ depending only on the dimension so that, for
every origin-symmetric convex body $K$ in $\R^n$ and every non-negative integrable function $f$ on $K$
\begin{equation}\label{main-problem}
\int_K f(x) dx \ \le\ T_n \ |K|^{1/n} \max_{\xi \in S^{n-1}} \int_{K\cap \xi^\bot} f(x) dx\ ?
\end{equation}
In other words, is it true that the sup-norm of the Radon transform of any probability density on a convex body of volume one is bounded
from below by a constant depending only on the dimension? The case where $f\equiv 1$ corresponds to the slicing problem of Bourgain.

It was proved in \cite{Koldobsky-2014, Koldobsky-2014b} that the answer to this question is affirmative
with $T_n\le O(\sqrt{n}).$ A different proof, based on the Blaschke-Petkantschin formula (see \cite{Schneider-Weil-book}) was given in \cite{Chasapis-Giannopoulos-Liakopoulos-2017}. Inequality (\ref{main-problem}) holds true with an absolute constant in place of $T_n$ for
intersection bodies, unconditional convex bodies and duals of convex bodies with bounded volume ratio \cite{Koldobsky-2015b},
and for the unit balls of $n$-dimensional subspaces of $L_p,\ p>2,$ with $C=O(\sqrt{p}),$ \cite{Koldobsky-Pajor-2017} (note that the
unit balls of subspaces of $L_p$ with $0<p\le 2$ are intersection bodies).

These results follow from a more general inequality proved in \cite{Koldobsky-2015b} for any
origin-symmetric star body $K$ in $\R^n,$ and any integrable non-negative function $f$ on $K$,
\begin{equation} \label{slicing-ovr}
\int_K f(x) dx \le  2\ d_{\rm {ovr}}(K,{\mathcal{I}}_n)  |K|^{\frac 1n} \max_{\xi \in S^{n-1}} \int_{K\cap \xi^\bot} f(x) dx.
\end{equation}
Now assuming that $K$ is an origin-symmetric convex body, by John's theorem \cite{John-1948} we get $d_{\rm {ovr}}(K,{\mathcal{I}}_n)\le \sqrt{n}.$
Also, the distance is bounded by an absolute constant for unconditional convex bodies \cite{Koldobsky-2015b} and for the unit balls of subspaces
of $L_p,\ p>2$  \cite{EMilman-2006, Koldobsky-Pajor-2017}. Clearly, if $K$ is an intersection body, the distance is 1. The proof of the inequality
(\ref{slicing-ovr}) in \cite{Koldobsky-2015b} is based on a stability result for sections of star bodies. However, in Section \ref{stab} we present the proof of a more general result which implies (\ref{slicing-ovr}).

The estimate $T_n\le O(\sqrt{n})$ is optimal. Klartag and the second named author showed in \cite{Klartag-Koldobsky} that there exists an
origin-symmetric convex body $M$ in $\R^n$ and a probability density $f$ on $M$ so that
$$\int_{M\cap H}f \le c\frac{\sqrt{\log\log n}}{\sqrt{n}} |M|^{-1/n},$$
for every affine hyperplane $H$ in $\R^n,$ where $c$ is an absolute constant. The convex body $M$ which provides
the example is a Gluskin-type random polytope generated by properly scaled random vectors $\theta_j$ uniformly distributed on
the sphere $S^{n-1}$, while the density $f$ is the density of an appropriate convolution of the standard
Gaussian measure on ${\mathbb R}^n$ with the sum of the Dirac masses of the $\theta_j$'s, restricted on $M$.
The logarithmic term was later removed by Klartag and Livshyts \cite{Klartag-Livshyts-2020}, who added a
``random rounding" technique to the previous construction. So, finally $T_n\ge  c\ \sqrt{n}.$

Another estimate of this kind, involving also non-central sections, was proved in \cite{Bobkov-Klartag-Koldobsky-2018}. Namely, there exists
an absolute constant $C$ so that
\begin{equation}\label{slicing}
\int_K f(x)\,dx \, \le \, C \sqrt{p}\ d_{\rm ovr}(K,L_p^n)\ |K|^{1/n}
\sup_{H} \int_{K\cap H} f(x)\,dx,
\end{equation}
for any $p\ge 1$, $n\in \N$, and any origin-symmetric convex body $K$ in $\R^n,$ where
$d_{\rm {ovr}}(K,L_p^n)$ is the outer volume ratio distance from $K$ to the class $L_p^n$ of the unit
balls of $n$-dimensional subspaces of $L_p([0,1]),$ and the supremum is taken over all affine hyperplanes
$H$ in $\R^n.$

\subsection{A quotient inequality for sections of functions}\label{stab}

A general inequality which implies both (\ref{bp1}) and (\ref{slicing-ovr}) was proved in \cite{Giannopoulos-Koldobsky-Zvavitch-2022}.
Let $K$ and $L$ be star bodies in $\R^n,$
and let $f,g$ be non-negative continuous functions on $K$ and $L$, respectively, with $\|g\|_\infty=g(0)=1.$
Then
$$
\frac{\int_Kf}{\left(\int_L g\right)^{\frac{n-1}n}|K|^{\frac 1n}}  \le \frac n{n-1} d_{\rm ovr}(K,{\mathcal{I}}_n)
\max_{\xi\in S^{n-1}} \frac{\int_{K\cap \xi^\bot} f}{\int_{L\cap \xi^\bot} g}.
$$
In fact, for any integer $0<k<n$ we have that
\begin{equation}\label{quotient}
\frac{\int_Kf}{\left(\int_L g\right)^{\frac{n-k}n}|K|^{\frac kn}}  \le \frac n{n-k} \left(d_{\rm ovr}(K,{\mathcal{BP}}_k^n)\right)^k
\max_{H\in G_{n,n-k}} \frac{\int_{K\cap H} f}{\int_{L\cap H} g}.
\end{equation}
For the proof of \eqref{quotient} fix $\delta>0$ and let $D\in {\mathcal{BP}}_k^n$ be a body such that $K\subset D$ and
\begin{equation}\label{sect11}
|D|^{\frac 1n}\le (1+\delta)\ d_{\rm ovr}(K,{\mathcal{BP}}_k^n)\ |K|^{\frac 1n}.
\end{equation}
Write $\nu_D$ for the measure on $G_{n,n-k}$ corresponding to $L$ according to the definition (\ref{genint}).
Let $\e$ be such that
$$\int_{K\cap H} f\le \e \int_{L\cap H} g,\qquad \mbox{for all}\;\; H\in G_{n,n-k}.$$
By (\ref{meas-sect}), we have
$$\Rad_{n-k}\left(\int_0^{\|\cdot\|_K^{-1}} r^{n-k-1}f(r\ \cdot)\ dr \right)(H) \le
\e\ \Rad_{n-k}\left(\int_0^{\|\cdot\|_L^{-1}} r^{n-k-1}g(r\ \cdot)\ dr \right)(H)$$
for every $H\in G_{n,n-k}.$
Integrating both sides of the latter inequality with respect to $\nu_D$ and using the definition (\ref{genint}), we get
\begin{align}\label{integration}
&\int_{S^{n-1}} \|x\|_D^{-k} \left(\int_0^{\|x\|_K^{-1}} r^{n-k-1}f(rx)\ dr \right)dx\\
\nonumber &\hspace*{1cm}\le \e \int_{S^{n-1}} \|x\|_D^{-k} \left(\int_0^{\|x\|_L^{-1}} r^{n-k-1}g(rx)\ dr \right)dx,
\end{align}
which is equivalent to
\begin{equation}\label{sect12}
\int_K \|x\|_D^{-k}f(x)dx \le \e \int_L \|x\|_D^{-k}g(x)dx.
\end{equation}
Since $K\subset D,$ we have $1\ge \|x\|_K\ge \|x\|_D$ for every $x\in K.$ Therefore,
$$\int_K \|x\|_D^{-k}f(x) dx \ge \int_K \|x\|_K^{-k}f(x) dx \ge \int_K f.$$
On the other hand, we may apply
\cite[Lemma 2.1]{VMilman-Pajor-1989}. Indeed,  recall that $g(0)=\|g\|_{\infty}=1$ and the same is true for the function $\chi_L(x)g(x)$, moreover
 the proof of \cite[Lemma 2.1]{VMilman-Pajor-1989}, also works under assumption that the body is star-shaped. Thus, we get
$$
\left(\frac{\int_{L}\|x\|_D^{-k} {g}(x) dx}{ \int_{D}\|x\|_D^{-k} dx} \right)^{1/(n-k)}  \le \left(\frac{\int_{L}{g}(x) dx}{\int_{D} dx} \right)^{1/n}.
$$
Since $\int_D\|x\|_D^{-k} dx =\frac{n}{n-k} |D|,$ we can estimate the right-hand side of
(\ref{sect12}) by
$$\int_L \|x\|_D^{-k} g(x) dx \le\e \frac n{n-k}  \left(\int_Lg\right)^{\frac {n-k}n} |D|^{\frac kn}.$$
Applying (\ref{sect11}) and sending $\delta$ to zero, we see that the latter inequality in conjunction with (\ref{sect12})
implies
$$\int_K f \le \e\ \frac n{n-k} \left(d_{\rm ovr}(K,{\mathcal{BP}}_k^n)\right)^k \left(\int_Lg\right)^{\frac {n-k}n} |K|^{\frac kn}.$$
Now put
$\e = \max\limits_{H\in G_{n,n-k}} \frac{\int_{K\cap H} f}{\int_{L\cap H} g}$ and the result follows.

Note the following immediate consequences of \eqref{quotient}. If we add the assumption that
$$
\int_{K\cap H}f\le \int_{L\cap H} g
$$
for all $H\in G_{n,n-k}$ then we get the generalization of \eqref{bp1}:
$$
\int_K f \le \frac n{n-k} \left(d_{\rm ovr}(K,BP_k^n)\right)^k |K|^{\frac kn}\left(\int_Lg\right)^{\frac{n-k}n}.
$$
If we choose $L=B_2^n$ and $g\equiv 1$ then we obtain the slicing inequality
$$\int_Kf\le \frac n{n-k} \left(d_{\rm ovr}(K,BP_k^n)\right)^k |K|^{\frac kn}\max_H \int_{K\cap H} f,$$
which generalizes \eqref{slicing-ovr}.
Choosing $K=L$ and $g\equiv 1$ we obtain another variant of the slicing inequality for functions:
$$\frac{\int_Kf}{|K|}\le \frac n{n-k} \left(d_{\rm ovr}(K,BP_k^n)\right)^k \max_H \frac{\int_{K\cap H}f}{|K\cap H|}.$$

\section{Projections of convex bodies}\label{projections}

\subsection{The Shephard's problem}

Shephard's problem \cite{Shephard-1964} is ``dual" to the Busemann-Petty problem: let $K$ and $L$ be two origin-symmetric
convex bodies in ${\mathbb R}^n$ and suppose that
\begin{equation}\label{eq:intro-4}|P_{\xi^{\perp } }(K)|\leq   |P_{\xi^{\perp} }(L)|\end{equation}
for every $\xi \in S^{n-1}$. Does it follow that $|K|\leq |L|$?

The answer is affirmative if $n=2$, but shortly after it was posed, Shephard's question was answered in the negative for all $n\geq 3$.
This was done independently by Petty in \cite{Petty-1967} who gave an explicit counterexample in ${\mathbb R}^3$, and by Schneider
in \cite{Schneider-1967} for all $n\geq 3$. In particular, Schneider
in \cite{Schneider-1967}   proved that the answer is affirmative if the body $L$ having projections of larger volume is a projection body; we refer to \cite{Koldobsky-Ryabogin-Zvavitch-2004} for Harmonic Analysis proofs of these facts.
After these counterexamples, one might try to relax the question, asking for the smallest constant
$C_n$ (or the order of growth of this constant $C_n$ as $n\to \infty $) for which: if $K,L$ are origin-symmetric convex bodies in
${\mathbb R}^n$ and $|P_{\xi^{\perp }}(K)|\leq   |P_{\xi^{\perp }}(L)|$ for all $\xi \in S^{n-1}$ then $|K|\leq   C_n|L|$.
Such a constant $C_n$ does exist, and a simple argument, based on John's theorem, shows that $C_n\leq   c\sqrt{n}$,
where $c>0$ is an absolute constant. On the other hand, K. Ball has proved in \cite{Ball-1991} that this simple
estimate is optimal: one has $C_n\approx\sqrt{n}$.

The lower dimensional Shephard problem is the following question. Let $1\leq   k\leq   n-1$
and let $S_{n,k}$ be the smallest constant $S >0$ with the following property: For
every pair of convex bodies $K$ and $L$ in ${\mathbb R}^n$ that satisfy $|P_F(K)|\leq   |P_F(L)|$
for all $F\in G_{n,n-k}$, one has $|K|^{\frac{1}{n}}\leq   S\,|L|^{\frac{1}{n}}$.
Is it true that there exists an absolute constant $C>0$ such that $S_{n,k}\leq C$ for all $n$ and $k$?

Goodey and Zhang \cite{Goodey-Zhang-1998} proved that $S_{n,k}>1$ if $n-k>1.$ General estimates for $S_{n,k}$ are
provided in \cite{Giannopoulos-Koldobsky-2017}: If $K$ and $L$ are two convex bodies in ${\mathbb R}^n$ such that
$|P_F(K)|\leq   |P_F(L)|,$ for every $F\in G_{n,n-k}$ then
\begin{equation*}
|K|^{\frac{1}{n}}\leq   c_1\sqrt{\frac{n}{n-k}}\log\left (\frac{en}{n-k}\right )\,|L|^{\frac{1}{n}},
\end{equation*}
where $c_1>0$ is an absolute constant. It follows that $S_{n,k}$ is bounded by an absolute constant if $\frac{n}{n-k}$ is bounded.
Also, under the same assumptions and using results from \cite{Paouris-Pivovarov-2013} one can prove a general estimate which
is logarithmic in $n$ and valid for all $k$:
\begin{equation*}|K|^{\frac{1}{n}}\leq   \frac{c_1\,\min w(\tilde{L})}{\sqrt{n}}\,|L|^{\frac{1}{n}}\leq   c_2(\log n)|L|^{\frac{1}{n}},\end{equation*}
where $c_1,c_2>0$ are absolute constants and the minimum is over all linear
images $\tilde{L}$ of $L$ that have volume $1$. The second inequality follows from the fact that if $\tilde{L}$ is a convex body of volume $1$ in ${\mathbb R}^n$ which is in the minimal mean width position (i.e. $w(\tilde{L})\leq   w(T(\tilde{L}))$ for all $T\in SL(n)$), then $w(\tilde{L})\leq   c\sqrt{n}(\log n)$ for some absolute constant $c>0$. This is a consequence of well-known results of Lewis, Figiel and Tomczak-Jaegermann, Pisier (see \cite[Chapter 6]{AGA-book} for a complete discussion).

An extension of Shephard's problem to the case of more general measures was first considered by Livshyts \cite{Liv-2016} who studied the case of $p$-concave and $1/p$-homogeneous measures.  Kryvonos and Langharst \cite{ Kryvonos-Langharst-2021} further extended the results from \cite{Liv-2016}, as well as, studied the isomorphic case of the question.

\subsection{A variant of the Busemann-Petty and Shephard problem}

A variant of the Busemann-Petty and Shephard problems was proposed by V.~Milman: Assume that $K$ and $L$ are origin-symmetric convex bodies in ${\mathbb R}^n$ and satisfy
\begin{equation}\label{milman}
|P_{\xi^{\perp }}(K)|\leq   |L\cap\xi^{\perp }|
\end{equation}
for all $\xi\in S^{n-1}$. Does it follow that $|K|\leq   |L|$? An affirmative answer to this
question was given by the first two authors in \cite{Giannopoulos-Koldobsky-2017}.
Also the lower dimensional analogue of the problem has an affirmative
answer, and moreover, one can drop the symmetry assumptions and even the assumption of convexity for $L$.
More precisely, if $K$ is a convex body in ${\mathbb R}^n$ and $L$ is a compact subset of ${\mathbb R}^n$ such that,
for some $1\leq   k\leq   n-1$,
\begin{equation*}|P_F(K)|\leq   |L\cap F|\end{equation*}
for all $F\in G_{n,n-k}$, then
\begin{equation*}|K|\leq   |L|.\end{equation*}
The proof exploits the Busemann-Straus/Grinberg inequality (see \cite{Busemann-Straus-1960}, \cite{Grinberg-1990})
\begin{equation}\label{eq:grin-1}\int_{G_{n,k}}|K\cap E|^nd\nu_{n,k}(E)\leq   \frac{\omega_k^n}{\omega_n^k}\,
|K|^k.\end{equation}
which is true for every bounded Borel set $K$ in ${\mathbb R}^n$ and $1\leq   k \leq   n-1$, and the classical
Alexandrov's inequalities in the following form: If $K$ is a convex body in ${\mathbb R}^n,$ then the sequence
\begin{equation}\label{eq:alexandrov-1}Q_k(K)=\left (\frac{1}{\omega_k}\int_{G_{n,k}}|P_F(K)|\,d\nu_{n,k}(F)\right)^{1/k}\end{equation}
is decreasing in $k$. This is a consequence of the Alexandrov-Fenchel inequality (see the books of Burago and Zalgaller \cite{Burago-Zalgaller-book}
and Schneider \cite{Schneider-book}). In particular, for every $1\leq   k\leq   n-1,$ we have
\begin{equation}\label{eq:alexandrov-2}
\left (\frac{|K|}{\omega_n}\right )^{\frac{1}{n}}\leq   \left  (\frac{1}{\omega_{n-k}}\int_{G_{n,n-k}}|P_F(K)|\,d\nu_{n,n-k}(F)\right )^{\frac{1}{n-k}}\leq   w(K).
\end{equation}
With these tools one proceeds as follows: Let $K$ be a convex body in ${\mathbb R}^n$ and $L$ be a compact subset of ${\mathbb R}^n$.
Assume that for some $1\leq   k\leq   n-1$ we have $|P_F(K)|\leq   |L\cap F|$ for all $F\in G_{n,n-k}$. From \eqref{eq:alexandrov-2} we get
\begin{equation*}\left (\frac{|K|}{\omega_n}\right )^{\frac{n-k}{n}}\leq   \frac{1}{\omega_{n-k}}\int_{G_{n,n-k}}|P_F(K)|\,d\nu_{n,n-k}(F).\end{equation*}
Our assumption, H\"{o}lder's inequality and the Busemann-Straus/Grinberg inequality give
\begin{align}\frac{1}{\omega_{n-k}}\int_{G_{n,n-k}}|P_F(K)|\,d\nu_{n,n-k}(F)
&\leq   \frac{1}{\omega_{n-k}}\int_{G_{n,n-k}}|L\cap F|\,d\nu_{n,n-k}(F)\\
&\hspace*{-1cm}\leq \frac{1}{\omega_{n-k}}\left (\int_{G_{n,n-k}}|L\cap F|^n\,d\nu_{n,n-k}(F)\right )^{1/n}\nonumber \\
&\hspace*{-1cm}\leq    \frac{1}{\omega_{n-k}}\,\frac{\omega_{n-k}}{\omega_n^{\frac{n-k}{n}}}|L|^{\frac{n-k}{n}}= \left (\frac{|L|}{\omega_n}\right )^{\frac{n-k}{n}}.\nonumber \label{eq:gk}
\end{align}
Therefore, $|K|\leq   |L|$.
\medbreak
Hosle \cite{Hosle-21} proved that if the condition (\ref{milman}) is reversed, then $|L|\le \sqrt{n}|K|.$ He also suggests a way to generalize the solution of Milman's problem to log-concave measures, $p$-concave measures and measures with $1/p$-homogeneous densities in place of volume.

\subsection{Surface area of projections}

The second named author obtained a number of ``slicing type" inequalities about the surface area of
hyperplane projections of projection bodies.
In \cite{Koldobsky-2013} he proved that if $Z$ is a projection body in ${\mathbb R}^n,$ then
\begin{equation}\label{eq:surface-K1}
|Z|^{\frac{1}{n}}\,\min_{\xi\in S^{n-1}}S(P_{\xi^{\perp }}(Z))\leq   b_nS(Z),
\end{equation}
where $S(A)$ denotes the surface area of $A$, and
\begin{equation*}b_n=(n-1)\omega_{n-1}/(n\omega_n^{\frac{n-1}{n}})\approx 1.\end{equation*}
This inequality is sharp; there is equality if $Z=B_2^n$.
Conversely, in \cite{Koldobsky-2015} he proved that if $Z$ is a projection body in ${\mathbb R}^n$ which is a dilate of a body in isotropic position (see, for example,  \cite{BGVV-book} ), then
\begin{equation}\label{eq:surface-K2}|Z|^{\frac{1}{n}}\,\max_{\xi\in S^{n-1}}S(P_{\xi^{\perp }}(Z))\geq c (\log n)^{-2}S(Z),\end{equation}
where $c>0$ is an absolute constant.

Similar inequalities for the surface area of hyperplane projections of an arbitrary convex body $K$ in ${\mathbb R}^n$
were studied in \cite{Giannopoulos-Koldobsky-Valettas-2018}. In what follows, we denote by $\partial_K$ the minimal surface
area parameter of $K$, which is the quantity
\begin{equation}\label{eq:partial}\partial_K:=\min\left\{ S(T(K))/|T(K)|^{\frac{n-1}{n}}:T\in GL(n)\right\}.\end{equation}
From the isoperimetric inequality and K.~Ball's reverse isoperimetric inequality \cite{Ball-1991c} it is known that
$c_1\sqrt{n}\leq   \partial_K\leq   c_2n$ for every convex body $K$ in ${\mathbb R}^n$, where $c_1,c_2>0$
are absolute constants. It was proved in \cite{Giannopoulos-Koldobsky-Valettas-2018}
that there exists an absolute constant $c_1>0$ such that, for every convex body $K$ in ${\mathbb R}^n$,
\begin{equation}\label{eq:intro-3}|K|^{\frac{1}{n}}\,\min_{\xi\in S^{n-1}}S(P_{\xi^{\perp }}(K))\leq  \frac{2b_n\partial_K}{n\omega_n^{\frac{1}{n}}}\,S(K)\leq  \frac{c_1\partial_K}{\sqrt{n}}\,S(K).\end{equation}
This inequality, which generalizes \eqref{eq:surface-K1}, is sharp e.g. for the Euclidean unit ball.
Since $c_1\partial_K/\sqrt{n}\leq   c\sqrt{n}$ for every convex
body $K$ in ${\mathbb R}^n$, we get the general upper bound
\begin{equation}|K|^{\frac{1}{n}}\,\min_{\xi\in S^{n-1}}S(P_{\xi^{\perp }}(K))\leq   c\sqrt{n}\,S(K).\end{equation}
The proof is based on the next result from
\cite{Giannopoulos-Hartzoulaki-Paouris-2002}. If $K$ is a convex body in ${\mathbb R}^n,$ then
\begin{equation}\label{eq:GHP}\frac{S(P_{\xi^{\perp }}(K))}{|P_{\xi^{\perp }}(K)|}\leq   \frac{2(n-1)}{n}\frac{S(K)}{|K|}\end{equation}
for every $\xi\in S^{n-1}$. It follows that
\begin{equation*}|K|\,\min_{\xi\in S^{n-1}}S(P_{\xi^{\perp }}(K))\leq   \frac{2(n-1)}{n}\,S(K)\,\min_{\xi\in S^{n-1}}|P_{\xi^{\perp }}(K)|.\end{equation*}
Next, we observe that
\begin{equation*}\min_{\xi\in S^{n-1}}|P_{\xi^{\perp }}(K)|=\min_{\xi\in S^{n-1}}h_{\Pi K}(\xi )=r(\Pi K),\end{equation*}
where $r(A)$ is the inradius of $A$, i.e. the largest $r>0$ such that $rB_2^n\subseteq A$. We write
\begin{equation*}r(\Pi K)\leq   \left (\frac{|\Pi K|}{\omega_n}\right )^{\frac{1}{n}}\leq   \frac{\omega_{n-1}\partial_K}{n\omega_n}\,|K|^{\frac{n-1}{n}},\end{equation*}
where the upper estimate for $|\Pi K|$ is observed in \cite{Giannopoulos-Papadimitrakis-1999}. Combining the above we get
\begin{equation*}|K|\,\min_{\xi\in S^{n-1}}S(P_{\xi^{\perp }}(K))\leq
\frac{2(n-1)\omega_{n-1}\partial_K}{n^2\omega_n}\,S(K)\,|K|^{\frac{n-1}{n}}.\end{equation*}
Inequality \eqref{eq:surface-K2} can be also generalized, starting with the next fact: If $K$ is a convex body in ${\mathbb R}^n$ then
\begin{equation*}\int_{S^{n-1}}S(P_{\xi^{\perp }}(K))\,d\sigma (\xi )\geq c\,S(K)^{\frac{n-2}{n-1}}.\end{equation*}
This inequality implies that if $K$ is in some of the ``classical positions" (we refer to \cite{AGA-book, VMilman-Pajor-1989} for discussion on those ``classical positions," including, minimal surface area, minimal mean width, isotropic,
John or L\"{o}wner position) then
\begin{equation}\label{eq:surface-GKV}|K|^{\frac{1}{n}}\,\int_{S^{n-1}}S(P_{\xi^{\perp }}(K))\,d\sigma (\xi )\geq c\,S(K),\end{equation}
where $c>0$ is an absolute constant. In particular, we get that if $K$ is a convex body in ${\mathbb R}^n$, which is in any of the ``classical positions"
then
\begin{equation*}|K|^{\frac{1}{n}}\,\max_{\xi\in S^{n-1}}S(P_{\xi^{\perp }}(K))\geq c\,S(K).\end{equation*}
Recall that a $(\log n)^2$-term appeared in \eqref{eq:surface-K2}. The estimate in \eqref{eq:surface-GKV} is stronger and,
for bounds of this type, there is no need to assume that $K$ is a projection body.
In fact, the estimate continues to hold as long as
\begin{equation*}S(K)^{\frac{1}{n-1}}\leq   c|K|^{\frac{1}{n}}\end{equation*}
for an absolute constant $c>0$. This is a mild condition which is satisfied not only by the classical positions but also
by all {\it reasonable} positions of $K$.

The same questions may be studied for the quermassintegrals $V_{n-k}(K)=V((K, n-k), (B_2^n,k))$ of a convex body $K$ and the
corresponding quermassintegrals of its hyperplane projections. The proofs employ the same tools as in the surface area case.
The main additional ingredient is a generalization of \eqref{eq:GHP} (proved in \cite{Fradelizi-Giannopoulos-Meyer-2003})
to subspaces of arbitrary dimension and quermassintegrals of any order: if $K$ is a
convex body in ${\mathbb R}^n$ and $0\leq   p\leq   k\leq   n$, then, for every $F\in G_{n,k}$,
\begin{equation}\label{mixedFGM}\frac{V_{n-p}(K)}{|K|}\geq\frac{1}{\binom{n-k+p}{n-k}}\frac{V_{k-p}(P_F(K))}{|P_F(K)|}.\end{equation}
This inequality allows one to obtain further generalizations;
one can compare the surface area of a convex body $K$ to the minimal, average or maximal surface area of its lower
dimensional projections $P_F(K)$, $F\in G_{n,k}$, for any given $1\leq   k\leq   n-1$. We refer to \cite{Fradelizi-Madiman-Meyer-Zvavitch2022, Fradelizi-Madiman-Zvavitch2022} to further study of inequalities related to (\ref{mixedFGM}).

\section{Surface area}\label{surface-area}

The question whether it is possible to have a version of the slicing inequality for the surface area instead of volume
has been formulated as follows: is it true that there exists a constant $\alpha_n$ depending (or not) on the dimension $n$ so that
\begin{equation}\label{eq:slicing-surface-1}S(K)\leq  \alpha_n|K|^{\frac{1}{n}}\max_{\xi\in S^{n-1}}S(K\cap\xi^{\perp })\end{equation}
for every origin-symmetric convex body $K$ in ${\mathbb R}^n$? A lower dimensional slicing problem may be also formulated; for any
$2\leq   k\leq   n-1$ one may ask for a constant $\alpha_{n,k}$ such that
\begin{equation}\label{eq:slicing-surface-lower}S(K)\leq  \alpha_{n,k}^k|K|^{\frac{k}{n}}\max_{H\in G_{n,n-k}}S(K\cap H)\end{equation}
for every origin-symmetric convex body $K$ in ${\mathbb R}^n$. Moreover, one may replace surface area by any other quermassintegral
and pose the corresponding question.

A negative answer was given in \cite{Brazitikos-Liakopoulos-2022}. For any $n\geq 2$ and any $\alpha>0$ one may find an
origin-symmetric convex body $K$ in ${\mathbb R}^n$ such that
\begin{equation*}S(K)>\alpha |K|^{\frac{1}{n}}\max\limits_{\xi\in S^{n-1}}S(P_{\xi^{\perp }}(K))\geq\alpha |K|^{\frac{1}{n}}\max\limits_{\xi\in S^{n-1}}S(K\cap\xi^{\perp }).\end{equation*}
In fact, it is shown that one may construct an origin-symmetric ellipsoid ${\mathcal E}$ such that
$$S({\mathcal E})>\alpha |{\mathcal E}|^{\frac{1}{n}}\max_{\xi\in S^{n-1}}S(P_{\xi^{\perp}}({\mathcal E})).$$
In order to do this, for a given ellipsoid ${\mathcal E}$ in ${\mathbb R}^n$ one needs to know the $(n-1)$-dimensional
section of ${\mathcal E}$ that has the largest surface area. This is a natural question of independent interest.
It is shown in \cite{Brazitikos-Liakopoulos-2022} that if ${\mathcal E}$ is an origin-symmetric ellipsoid in ${\mathbb R}^n$, and if
$a_1\leq a_2\leq\cdots \leq a_n$ are the lengths and $e_1,e_2,\ldots ,e_n$ are the corresponding directions of its semi-axes, then
\begin{equation}\label{eq:max-ellipsoid}S({\mathcal E}\cap \xi^{\perp })\leq S(P_{\xi^{\perp}}({\mathcal E})) \leq S({\mathcal E}\cap e_1^{\perp })\end{equation}
for every $\xi\in S^{n-1}$. This information is then combined with a formula of Rivin \cite{Rivin-2007}
for the surface area of an ellipsoid: If ${\mathcal E}$ is an ellipsoid as above then
\begin{equation}\label{eq:rivin}S({\mathcal E})=n\,|{\mathcal E}|\,\int_{S^{n-1}}\left (\sum_{i=1}^n\frac{\xi_i^2}{a_i^2}\right)^{1/2}d\sigma (\xi ).\end{equation}
Assume that there exists a constant $\alpha_n>0$ such that we have the following inequality for ellipsoids:
\begin{equation}\label{eq:max}S({\mathcal E})\leq \alpha_n|{\mathcal E}|^{1/n}\max_{\xi\in S^{n-1}}S({\mathcal E}\cap\xi^{\perp }).\end{equation}
Then we have
$$\max_{\xi\in S^{n-1}}S({\mathcal E}\cap\xi^{\perp })=S({\mathcal E}\cap e_1^{\perp })=(n-1)\,|{\mathcal E}\cap e_1^{\perp }|\,\int_{S^{n-2}}\Big (\sum_{i=2}^n\frac{\xi_i^2}{a_i^2}\Big)^{1/2}d\sigma (\xi )$$
and assuming, as we may, that $\prod_{i=1}^na_i=1$ we can rewrite \eqref{eq:max} as
\begin{equation*}n\omega_n\cdot\frac{1}{d_n}{\mathbb E}\Big[\Big (\sum_{i=1}^n\frac{g_i^2}{a_i^2}\Big)^{1/2}\Big]
\leq \alpha_n\omega_n^{1/n}\cdot (n-1)\omega_{n-1}\frac{1}{a_1}\cdot \frac{1}{d_{n-1}}{\mathbb E}\Big[\Big (\sum_{i=2}^n\frac{g_i^2}{a_i^2}\Big)^{1/2}\Big],\end{equation*}
where $d_n\sim\sqrt{n}$. After some calculations we see that
\begin{equation*}\alpha_n \geq c\,\left (\frac{1+\sum_{i=2}^n\frac{a_1^2}{a_i^2}}{\sum_{i=2}^n\frac{1}{a_i^2}}\right )^{1/2}
\end{equation*}
for some absolute constant $c>0$. Choosing $a_2=\cdots =a_n=r$ and $a_1=r^{-(n-1)}$ we see that
$$\left (\frac{1+\sum_{i=2}^n\frac{a_1^2}{a_i^2}}{\sum_{i=2}^n\frac{1}{a_i^2}}\right )^{1/2}=
\left (\frac{1+\frac{n-1}{r^{2n}}}{\frac{n-1}{r^2}}\right )^{1/2}=\left (\frac{1}{r^{2n-2}}+\frac{r^2}{n-1}\right )^{1/2}\to\infty $$
as $r\to\infty $. So, we arrive at a contradiction.

In fact, one can prove a more general analogue of \eqref{eq:max-ellipsoid}; for any $k-$dimensional subspace $H$ and any $0\leq j\leq k-1$ we have that
$$W_j({\mathcal E}\cap F_k)\leq W_j({\mathcal E}\cap H)\leq W_j({\mathcal E}\cap E_k)$$
and
$$W_j(P_{F_k}({\mathcal E}))\leq W_j(P_H({\mathcal E}))\leq W_j(P_{E_k}({\mathcal E})),$$
where $F_k=\mathrm{span}\{e_1,\ldots ,e_k\}$, $E_k=\mathrm{span}\{e_{n-k+1},\ldots, e_n\}$ and $W_j(A)=V((A,k-j),(B_2^k,k-j))$
is the $j$-th quermassintegral of a $k$-dimensional convex body $A$.
These results generalize a known fact for the maximal and minimal volume of $k$-dimensional sections and projections of ellipsoids.
As a consequence one can obtain a more general negative result about all the quermassintegrals
of sections and projections of convex bodies.

In \cite{Brazitikos-Liakopoulos-2022} some positive results are stated for variants of this question, which were strengthened
in \cite{Liakopoulos-preprint}. An example is the next inequality: If $K$ is an origin-symmetric convex body in ${\mathbb R}^n$
then for any $0\leq  j\leq  n-k-1\leq  n-1$ we have
$$\alpha_{n,k,j}\binom{n}{k}^{-1}\frac{W_j(K)}{|K|}\leq \int_{G_{n,n-k}}\frac{W_j(K\cap F)}{|K\cap F|}d\nu_{n,n-k}(F)
\leq  \alpha_{n,k,j}\binom{n-j}{k}\frac{W_j(K)}{|K|}$$
where $\alpha_{n,k,j}$ is a constant depending only on $n,k,j$.

The analogue of the Busemann-Petty problem for surface area was studied by Koldobsky and K\"{o}nig in \cite{KK-2019}:
If $K$ and $L$ are two convex bodies in ${\mathbb R}^n$ such that
$S(K\cap\xi^{\perp })\leq S(L\cap\xi^{\perp })$ for all $\xi\in S^{n-1}$ does it then follow that $S(K)\leq S(L)$? Answering a
question of Pelczynski, they prove that the central $(n-1)$-dimensional section of the cube $B_{\infty }^n=[-1,1]^n$ that has
maximal surface area is the one that corresponds to the unit vector $\xi_0=\frac{1}{\sqrt{2}}(1,1,0,\ldots ,0)$ (exactly as in
the case of volume) i.e.
$$\max_{\xi\in S^{n-1}}S(B_{\infty }^n\cap\xi^{\perp })=S(B_{\infty }^n\cap\xi_0^{\perp })=2((n-2)\sqrt{2}+1).$$
Comparing with a ball of suitable radius one gets that the answer to the Busemann-Petty problem for surface area is negative
in dimensions $n\geq 14$.

It is natural to ask whether an isomorphic version of the problem has an affirmative answer. Assuming that there is a
constant $\gamma_n$ such that if $K$ and $L$ are origin-symmetric convex bodies in $\mathbb{R}^n$ that satisfy
$$S(K\cap\xi^{\perp})\leq  S(L\cap\xi^{\perp})$$ for all $\xi\in S^{n-1}$, then $S(K)\leq  \gamma_n S(L)$,
one can see that there is some constant $c(n)$ such that
\begin{equation}\label{slicing-15}
S(K)\leq  c(n)S(K)^{\frac{1}{n-1}}\max_{\xi\in S^{n-1}} S(K\cap\xi^{\perp})
\end{equation}
for every convex body $K$ in ${\mathbb R}^n$. It was proved in \cite{Brazitikos-Liakopoulos-2022}
that an inequality of this type holds true in general. If $K$ is a convex body in ${\mathbb R}^n$ then
$$S(K)\leq  A_nS(K)^{\frac{1}{n-1}}\max\limits_{\xi\in S^{n-1}}S(K\cap\xi^{\perp })$$
where $A_n>0$ is a constant depending only on $n$. The result is first proved for an arbitrary ellipsoid
and then it is extended to any convex body, using John's theorem. In \cite{Liakopoulos-preprint} a direct proof of a more general
result is given, showing that an inequality as \eqref{slicing-15} holds for any $k$ and $j$, where $k$ is the codimension
of the subspaces and $j$ is the order of the quermassintegral that we consider:
Let $K$ be an origin-symmetric convex body in $\mathbb{R}^n$. For every $0\leq  j\leq  n-k-1\leq  n-1$ we have that
\begin{equation*}
W_j(K)^{n-k-j}\leq  \alpha_{n,k,j}\max_{F\in G_{n,n-k}}W_j(K\cap F)^{n-j},
\end{equation*}
where $\alpha_{n,k,j}>0$ is a constant depending only on $n,k$ and $j$.
The proof of this inequality exploits the Blaschke-Petkantschin formula
and some integral-geometric results of Dann, Paouris and Pivovarov
from \cite{Dann-Paouris-Pivovarov-2016}.

\section{Volume difference inequalities}\label{vol-diff}

Volume difference inequalities estimate the error in computations of volume of a body out of the areas of its sections and projections.
Starting with the case of sections, let $\gamma_{n,k}$ be the smallest constant $\gamma>0$ with the property that
\begin{equation} \label{main-prob1}
|K|^{\frac {n-k}n}-|L|^{\frac {n-k}n}\leq \gamma^k
\max_{F\in G_{n,n-k}} \big (|K\cap F|-|L\cap F|\big )
\end{equation}
for all $1\leq k <n$ and all origin-symmetric convex bodies $K$ and $L$ in $\R^n$ such that $L\subset K.$
The question is whether there exists an absolute constant $C$ so that $\sup_{n,k} \gamma_{n,k}\leq C$.
Note that without extra assumptions on $K$ and $L,$ inequality \eqref{main-prob1} cannot hold with any $\gamma>0,$
because of the counterexamples to the Busemann-Petty problem. Note also that if we apply \eqref{main-prob1}
with $L=\beta B_2^n$ as $\beta\to 0$, we get the slicing inequality.

It was proved in \cite{Koldobsky-2011} for $k=1$, and in \cite{Koldobsky-Ma-2013} for $1<k<n$ that if
$K\in {\mathcal{BP}}_k^n$ and $L$ is any origin-symmetric star body in $\R^n$ then \eqref{main-prob1} is true
in the form
\begin{equation} \label{initial-sect}
|K|^{\frac {n-k}n}-|L|^{\frac {n-k}n}\leq c_{n,k}^k
\max_{F\in G_{n,n-k}} \big(|K\cap F|-|L\cap F|\big),
\end{equation}
where $c_{n,k}^k=\omega_n^{\frac {n-k}n}/\omega_{n-k}$. One can check that $c_{n,k}\in (\frac 1{\sqrt{e}},1)$ for all $n,k.$

In \cite{Giannopoulos-Koldobsky-2018}, the inequality \eqref{initial-sect} was extended to arbitrary origin-symmetric star bodies.
Let $1\leq k <n,$ and let $K$ and $L$ be origin-symmetric star bodies in $\R^n$ such that
$L\subset K.$ Then
\begin{equation} \label{vdi-volume-ineq}
|K|^{\frac {n-k}n}-|L|^{\frac {n-k}n}\leq c_{n,k}^k d_{\rm ovr}^k(K,{\mathcal{BP}}_k^n)
\max_{F\in G_{n,n-k}} \big(|K\cap F|-|L\cap F|\big).
\end{equation}
The outer volume ratio distance was estimated in \cite{Koldobsky-Paouris-Zymonopoulou-2011}. If $K$ is an origin-symmetric convex
body in $\R^n,$ then
\begin{equation}\label{eq:ovr-10}d_{\rm ovr}(K,{\mathcal{BP}}_k^n)\leq c\sqrt{n/k}\,\big(\log (en/k)\big)^{\frac{3}{2}},\end{equation}
where $c>0$ is an absolute constant.  Therefore, \eqref{vdi-volume-ineq} provides
an affirmative answer to the question for sections of proportional dimension.

The volume difference inequality \eqref{vdi-volume-ineq} can be extended to arbitrary measures in place of volume,
as follows. Let $f$ be a bounded non-negative measurable function on $\R^n.$ Let $\mu$ be the measure
with density $f$ so that $\mu(B)=\int_B f$  for every Borel set $B$ in $\R^n.$ Also, for every $F\in G_{n,n-k}$ we write
$\mu(B\cap F)=\int_{B\cap F} f,$ where we integrate the restriction of $f$ to $F$ against the Lebesgue measure
on $F.$ For any $1\leq k <n$ and any pair of origin-symmetric star bodies $K$ and $L$ in $\R^n$ such that
$L\subset K,$ and any measure $\mu $ with even non-negative continuous density,
\begin{equation} \label{vdi-measure-ineq}
\mu(K)-\mu(L)\leq \frac n{n-k} c_{n,k}^k\ |K|^{\frac kn}\ d_{\rm ovr}^k(K,{\mathcal{BP}}_k^n)
\max_{F\in G_{n,n-k}} \big(\mu(K\cap F)-\mu(L\cap F)\big).
\end{equation}
In the opposite direction, for any measure in ${\mathbb R}^n$ with bounded density $g$,
\begin{equation}\big (\mu (K)-\mu (L)\big )^{\frac{n-k}{n}} \geq \frac{c_{n,k}^k}{\|g\|_{\infty }^{\frac{k}{n}}}
\left (\int_{G_{n,n-k}}\big (\mu (K\cap F)-\mu (L\cap F)\big )^{\frac{n}{n-k}}\,d\nu_{n,n-k}(F)\right )^{\frac{n-k}{n}}.\end{equation}
In particular,
\begin{equation}\big (\mu (K)-\mu (L)\big )^{\frac{n-k}{n}}
\geq c_{n,k}^k\frac{1}{\|g\|_{\infty }^{\frac{k}{n}}}\min_{F\in G_{n,n-k}}\big (\mu (K\cap F)-\mu (L\cap F)\big ).\end{equation}
There is also a result for projections. Let $\beta_n$ be the smallest constant $\beta>0$ satisfying
\begin{equation}\label{question-proj}
\beta \big (|L|^{\frac {n-1}n}-|K|^{\frac {n-1}n}\big ) \geq
\min_{\xi\in S^{n-1}} \,\big(|P_{\xi^\bot}(L)|-|P_{\xi^\bot}(K)|\big)
\end{equation}
for all origin-symmetric convex bodies $K,L$ in $\R^n$ whose curvature
functions $f_K$ and $f_L$ exist and satisfy
$f_K(\xi)\leq f_L(\xi)$ for all $\xi\in S^{n-1}$. Then, $\beta_n\simeq \sqrt{n},$
i.e. there exist absolute constants $a,b>0$ such that for all $n\in \N$
$$a\sqrt{n}\leq \beta_n \leq b\sqrt{n}.$$
Note that without an extra condition on $K$ and $L$, \eqref{question-proj}
cannot hold in general with any $\beta>0$, because of the counterexamples to the Shephard problem.

It was proved in \cite{Koldobsky-2011, Koldobsky-2013} that if $L$ is a projection body
and $K$ is an origin-symmetric convex body, then
\begin{equation} \label{initial-proj}
|L|^{\frac {n-1}n}-|K|^{\frac {n-1}n} \geq c_{n,1}
\min_{\xi\in S^{n-1}} \,\big(|P_{\xi^\bot}(L)|-|P_{\xi^\bot}(K)|\big).
\end{equation}
The condition $f_K\leq f_L$ is not needed for \eqref{initial-proj} because we assume
that $L$ is a projection body. This inequality is extended in \cite{Giannopoulos-Koldobsky-2018}
to arbitrary origin-symmetric convex bodies, as follows. Suppose that $K$ and $L$ are origin-symmetric convex
bodies in $\R^n,$ and their curvature functions exist and satisfy $f_K(\xi)\leq f_L(\xi)$ for all $\xi\in S^{n-1}.$
Then
\begin{equation} \label{vdi-proj}
d_{\rm {vr}}(L,\Pi)\,\big(|L|^{\frac {n-1}n}-|K|^{\frac {n-1}n}\big ) \geq c_{n,1}
\min_{\xi\in S^{n-1}} \,\big(|P_{\xi^\bot}(L)|-|P_{\xi^\bot}(K)|\big).
\end{equation}
Again by K.~Ball's volume ratio estimate, for any convex body $K$ in $\R^n,$
$d_{\rm vr}(K,\Pi)\leq \sqrt{n}.$ In fact, this distance
can be of the order $\sqrt{n},$ up to an absolute constant.

For more on volume difference inequalities, see \cite{KW}.

\section{Discrete versions}\label{discrete}

Let $\Z^n$ be the standard integer lattice in $\R^n$. Given an origin-symmetric convex body $K$, define $\#K={\rm card}(K\cap  \Z^n)$, the number of points of $\Z^n$ in $K$. The original  proof of the  Loomis-Whitney inequality \eqref{LW} is based on a discretization technique, i.e. to
consider a set $K$ which is decomposed into a union of equal disjoint cubes and restate the question in combinatorial language.
In particular, this combinatorial version implies (and is actually equivalent to) the following discrete variant
\begin{equation}\label{CLW}
\#K \le \left(\prod\limits_{i=1}^n \# (P_{e_i^\perp}(K)) \right)^{\frac{1}{n-1}}
\end{equation}
of the Loomis-Whitney inequality.

There are a number of very interesting tomographic questions related to the number of integer points in
the projections of a convex body; we refer to \cite{GGroZ, RYZh,  N-Zhang}.
We also note that in the recent years there are several attempts to translate questions and facts from classical convexity to more general
settings including discrete geometry. Here we will concentrate on inequalities concerning sections of a convex body.
The properties of sections of convex bodies with respect to the integer lattice were extensively studied in discrete tomography,
see \cite{GGroZ, GGr1, GGr2, GGro, RYZh},  where many interesting new properties were proved and a series of exciting open questions
were proposed. It is interesting to note that after translation many questions become quite non-trivial and counterintuitive,
and the answer may be quite different from the one in the continuous case. In addition, finding the relation between the geometry
of a convex set and the number of integer points contained in the set is always a non-trivial task. One can see this, for example,
from the history of Khinchin's flatness theorem (see \cite{Ba1, Ba2, BLPS, KL}).

Around  2013 the second named author asked if it is possible to provide a discrete analogue of the inequality \eqref{main-problem}:

\smallskip

\noindent {\it Question~1.} Does there exist a  constant $d_n$ such that
$$
\#K\leq d_n \max_{\xi\in S^{n-1}} \left( \#(K\cap \xi^\perp)\right)  |K|^{\frac{1}{n}},
$$
for all origin-symmetric convex bodies  $K\subset \R^n$ containing $n$ linearly independent lattice points?

We note here that we require that $K$ contains $n$ linearly independent lattice points in order to eliminate the
degenerate case of a body (for example, a box of the form $[-\delta,\delta]^{n-1}\times [-20, 20]$) whose maximal section contains all
lattice points in the body, but whose volume may be arbitrarily close to $0$ by considering sufficiently small $\delta>0$.

Thus, the methods applied to attack this question are quite different from the methods described in the previous sections and seem
to require use of tools from the geometry of numbers. Let us start with the simplest case and show that the constant $d_2$ exists,
i.e. it is independent from the origin-symmetric planar convex body $K$. We will use two classical statements (see for example, \cite{TV2}):

\smallskip

\noindent {\it Minkowski's First Theorem:} Let $K \subset \R^n$ be an origin-symmetric convex body such that $|K|\ge 2^n$.
Then $K$ contains at least one non-zero element of $\Z^n$.

\smallskip

\noindent {\it Pick's Theorem:} Let $P$ be an integral $2$-dimensional convex polygon. Then $|P|=I+\frac{1}{2}B-1$, where
$I$ is the number of lattice points in the interior of $P$, and $B$ is the number of lattice points on its boundary.
Here, a polygon is called integral if it can be described as the convex hull of its  lattice points.

\smallskip

Following \cite{AHZ-2017}, consider an origin-symmetric planar convex body $K$ and
let $s=\sqrt{|K|/4}$. By Minkowski's first theorem, there exists a non-zero vector $u\in \Z^2 \cap \frac{1}{s} K$.
Then $s u \in K$ and
$$\# \left(L_u \cap K\right) \geq 2 \lfloor s \rfloor + 1 
,$$ where $\lfloor s \rfloor$ is the integer part of $s$ and $L_u$ is the line containing $u$ and the origin.
Next, consider the convex hull $P$ of the lattice points inside $K$.  Since $P$ is an integral $2$-dimensional convex polygon,
by Pick's theorem we get that
$$
|P|=I+\frac{1}{2}B-1\geq \frac{I+B}{2}-\frac{1}{2},
$$
using that $I \ge 1$.  Thus,
$$
\# P=I+B   \leq 2 |P|+1\leq \frac{5}{2} |P|,$$
since the area of an origin-symmetric
integral convex polygon is at least $2$. It follows that
\begin{align*}
\# K  = \# P &\leq \frac{5}{2}|P| \leq \frac{5}{2}|K| = \frac{5}{2} \,(2\,s)\, |K|^{\frac{1}{2}}< 4\,( 2 \lfloor s \rfloor + 1)\,
|K|^{\frac{1}{2}}\\
&\leq  4 \max_{\xi\in S^{1}}  \#(K\cap \xi^\perp)\,|K|^{\frac{1}{2}}.
\end{align*}
Unfortunately, there seems to be no straightforward generalization of the above approach to higher dimensions.
This is partially due to the fact that the hyperplane sections of dimension higher than one are much harder to study,
but is also due to the lack of a direct analogue of Pick's formula. It is also essential to note that, in general,
the constant $d_n$ is dependent on $n$.  Indeed, for the cross-polytope $B_1^n= \{ x\in\R^n: \|x\|_1\leq 1\}$ we have
$\# B_1^n=2\,n+1$ and $\max\{\# (B_1^n\cap \xi^\perp):  \xi\in S^{n-1}\}=\#(B_1^n \cap e_1^\perp)=2n-1$, and $|B_1^n|^{1/n} \sim n^{-1}$.
Thus, $d_n$ must be greater than $cn$.  Using the special structure of unconditional bodies, it was proved in \cite{AHZ-2017}
that $d_n$ is of order $n$ for this class of convex bodies.

Another example which illustrates this situation is the classical Brunn's principle (see for example \cite{Gardner-book})
which tells us that for an origin-symmetric convex body $K$ one has
$$
|K\cap \xi^\perp| \ge  |K\cap (t\xi +\xi^\perp )|, \mbox{ for all } \xi \in S^{n-1}, t\in \R.
$$
The above statement is not true if volume is replaced by the number of integer points. Indeed, let $Q=[0,1]^{n-1} \subset \R^{n-1}$
and let $K$ be the convex hull of $Q+e_n$ and $-Q-e_n$. Then $\#(K\cap e_n^\perp)=1$ but $\#(K\cap (e_n+e_n^\perp ))=2^{n-1}.$
The following analogue of Brunn's concavity principle was proved in \cite{AHZ-2017,FH-2021}:
\begin{equation}\label{dbrunn}
\#(K\cap \xi^\perp) \ge  2^{1-n} \#(K\cap (t\xi +\xi^\perp )), \mbox{ for all } \xi \in S^{n-1}, t\in \R.
\end{equation}

To prove that the constant $d_n$ is bounded for a general origin-symmetric convex body, one may use the discrete analogue
of John's theorem \cite{BV-1992, TV1, TV2, BH-2019}, which gives an approximation of an origin-symmetric convex body by a symmetric
generalized arithmetic progression. This approach was used in \cite{AHZ-2017} to prove that $d_n \le O(n)^{7n/2}$.
The latter estimate can be slightly improved: using \cite{BH-2019} one can show that $d_n \le O(n)^{n},$ which is still far from optimal.
Finally, the following fact was proved in \cite{AHZ-2017}:
$$
\# K \leq O(1)^n n^{n-m} \max\left(\# (K\cap H)\right)\,\, |K|^{\frac{n-m}{n}},
$$
where the maximum is taken over all $m$-dimensional linear  subspaces $H \subset \R^n$. In particular, $d_n \le C^n$
for some large $C>0$. The proof of this fact is based on Minkowski's second theorem, the discrete Brunn's principle (\ref{dbrunn})
and the Bourgain-Milman inequality. We refer to \cite{AHZ-2017} for more details. Here we would like to present a probabilistic-harmonic
analysis approach to Question~1  which is due to Oded Regev \cite{Regev}.  We will show that there is a distribution of directions
$\xi$ for which $\#(K\cap \xi^\perp)$ is large enough, with a positive probability. For the convenience of the reader
we need to provide a few standard definitions and to prove some technical estimates. A lattice $\Lambda \subset \R^n$ is the set of all integer linear combinations of $n$ linearly independent vectors in $\R^n$. We notice that  $ \Lambda= T \Z^n$, where $T\in GL(n)$, $\det(T)\neq 0$ and
denote $\det (\Lambda)=\det(T)$.  We will also consider the dual lattice $\Lambda^*=\{x\in \R^n:  \langle x, y\rangle \in \Z, \mbox{ for all } y \in \Lambda\}$. Finally  for any $s>0$  consider the function $\rho_s(x): \R^n \to (0,1]$, defined by
$
\rho_s(x)=e^{-\pi |x|^2/s^2},
$
where $|x|$ is the Euclidean norm of $x$, and for a countable set $A\subset \R^n$ let
$$
\rho_s(A) =\sum\limits_{x \in A} e^{-\pi |x|^2/s^2}.
$$
We will use the Poisson summation formula (see for example \cite[Lemma~17.2]{Bar}). Consider a  function $f:\R^n \to \C$ and let
$$
\widehat{f}(y)=\int_{\R^n} e^{-2\pi i \langle x, y\rangle } f(x) dx \qquad (y\in\R^n)
$$
be the Fourier transform of $f$. Assume that $f$ and $\widehat{f}$ decay sufficiently fast, i.e. there exist positive constants
$C, \delta$ such  that $|f(x)|, |\widehat{f}(x)| \le \frac{C}{1+|x|^{n+\delta}}$ for all $x\in \R^n$ (notice that this condition
it trivially satisfied by $\rho_s$). Then,
\begin{equation}\label{poisson}
\sum\limits_{x \in \Lambda} f(x) =\frac{1}{\det(\Lambda)} \sum\limits_{y \in \Lambda^*} \widehat{f}(y).
\end{equation}
Using the fact that $\widehat{\rho}_s(y)=s^n\rho_{1/s}(y)$ together with (\ref{poisson}) we get
$$
\rho_s(\Lambda)=\sum\limits_{x \in A} \rho_s(x) =\frac{1}{\det(\Lambda)} \sum\limits_{y \in \Lambda^*} \widehat{\rho}_s(y) = \frac{s^n}{\det(\Lambda)}\rho_{1/s}(\Lambda^*).
$$
In particular, we see that
\begin{equation}\label{est1}
\rho_s(\Lambda) \geq {\det (\Lambda)^{-1}}\, s^n.
\end{equation} We may use again (\ref{poisson}) to show that $\rho_s(\Lambda +a) \le \rho_s(\Lambda)$ for all $a\in \R^n$.
Indeed, using $\widehat{\rho}_s(y)=s^n\rho_{1/s}(y) >0,$ we get
\begin{equation}\label{shift}
\rho_s(\Lambda +a)=\frac{1}{\det(\Lambda)}\!\! \sum\limits_{y \in \Lambda^*} \widehat{\rho}_s(y)  e^{-2\pi i \langle a, y\rangle } \!\!\le \!\!  \frac{1}{\det(\Lambda)}\!\! \sum\limits_{y \in \Lambda^*} \widehat{\rho}_s(y)  = \rho_s(\Lambda).
\end{equation}
Let $Z_{A,s}$ be a random vector taking values in a countable set $A$, such that
$$
\Pr(Z_{A,s}=x)=\frac{\rho_s(x)}{\rho_s(A)}.
$$
Using (\ref{est1}) we get
\begin{equation}\label{zero}
\Pr(Z_{\Lambda,s}=0)=\frac{1}{\rho_s(\Lambda)}\le \frac{\det(\Lambda)}{s^n}.
\end{equation}
Finally, if we pick a lattice $\Lambda \subset \R^n,$  $x\in \Lambda$ and $s>0$, we claim that
\begin{equation}\label{section}
\Pr(  \langle Z_{\Lambda^*,s} , x\rangle =0) \ge \Pr(Z_{\Z\setminus |x|,s} =0) =\rho_{s|x|}(\Z)^{-1} \ge c\min\{1, (s|x|)^{-1}\}.
\end{equation}
We may assume that $x\neq 0$ and note that by definition $\langle Z_{\Lambda^*,s} , x\rangle $ can take only integer values. Fix some $k\in \Z$ and
consider the set of points $y \in \Lambda^*$ such that $\langle x,y\rangle =k$. Note that if this set is empty then its $\rho_s$ mass is clearly zero.
Now consider the case where this set is not empty. The affine subspace $\{y \in \R^n: \langle x,y\rangle =k\}$ is a shift of $x^\perp$ at distance
$k/|x|$ in the direction of $x/|x|$. Notice that such a shift will not necessarily move points of $\Lambda^*$ into itself, and thus an
additional shift inside $\{y \in \R^n: \langle x,y\rangle =k\} \cap \Lambda^*$ may be required. Thus, using the product structure of $\rho_s$
and property (\ref{shift}) we get
$$
\rho_s(\{y \in \R^n: \langle x,y\rangle =k\} \cap \Lambda^*) \le \rho_s(k/|x|)\rho_s(\Lambda^*\cap x^\perp).
$$
Dividing the above inequality by $\rho_s(\Lambda^*)$ and summing up over $k \in \Z$ we get (\ref{section}).

The following theorem is due to O.~Regev \cite{Regev}: Consider an origin-symmetric convex body $K \subset \R^n$. Then,
$$
\#K\leq \max_{\xi\in S^{n-1}} \left( \#(K\cap \xi^\perp)\right) \max\{ 1, c n |K|^{\frac{1}{n-1}}\}.
$$
Note that this theorem improves the bound on $d_n$ from \cite{AHZ-2017}, when the volume of $|K|$ is smaller then $C^{n^3}$
and provides a polynomial bound on $d_n$ for bodies of volume smaller then $n^{cn^2}$.

For the proof, using John's theorem we see that there is linear transformation such that $\det(T)=1$ and $TK \subset n |TK|^{1/n} B_2^n$. Without loss of generality we will consider the body $TK$ instead of $K$ and the lattice $\Lambda=T \Z^n$ instead of $\Z^n$. We may also assume that $|K|\ge n^{-n}$,  otherwise $K\subset \delta B_2^n,$ where $0<\delta<1$ and using $\det(\Lambda)=1$ we get that $K\cap \Lambda$ is not full dimensional and the statement is trivial.

Now we will select the direction $\xi \in \Lambda^*\setminus\{0\}$, using a probabilistic approach. Let $\xi=Z_{\Lambda^*,s}$ where
$s \ge 1$ will be chosen later. Then we may apply (\ref{section}) to claim that for any fixed $x \in K \cap \Lambda$ we have
$$
\Pr( x \in \xi^\perp)\ge c \min\{1, (s|x|)^{-1}\} \ge  c \min\{1, (s n|K|^{1/n})^{-1}\} =\frac{c}{s n|K|^{1/n}}
$$
Then
$$
\Exp \left(\frac{\#\{x \in K\cap \Lambda: \langle x,\xi\rangle =0\}}{\#(K\cap \Lambda)}\right)\ge \frac{ c}{  s n|K|^{1/n}}. 
$$
The above inequality allows us to select $\xi$ for which $\#\{x \in K\cap \Lambda: \langle x,\xi\rangle =0\}$ is large.  To finish the proof, we need to
choose $s \ge 1 $ to guarantee that  $\xi$ can be selected not to be equal to $0$.  Using (\ref{zero}) we see that
$\Pr(\xi =0)\le s^{-n}.$  Thus we need to pick $s\ge 1$ such that $s^{-n} \le  c  (s n|K|^{1/n})^{-1}$, i.e. $s=C |K|^{\frac{1}{n(n-1)}} \ge 1$, where $C>0$ is large enough absolute constant. This completes the proof.

Another very interesting bound on the cardinality of lattice points in sections of convex bodies
is inspired by Meyer's inequality \eqref{eq:LW-2}. It was proposed by Gardner, Gronchi and Zong \cite{GGroZ} and proved by
Freyer and Henk in \cite{FH-2021}: For any origin-symmetric convex body $K\subset \R^n$ there exists
a basis $b_1, \dots, b_n$  of $\Z^n$ such that
$$
\left(\#K\right)^{\frac{n-1}{n}} \leq  O(n^22^n) \prod\limits_{i=1}^n  \left( \#(K\cap b_i^\perp)\right)^{1/n}
$$
and there are $t_1,\ldots ,t_n \in \Z^n$ such that
$$
\left(\#K\right)^{\frac{n-1}{n}} \leq O(n^2) \prod\limits_{i=1}^n  \left( \#(K\cap (t_i+ b_i^\perp))\right)^{1/n}.
$$
Thus we immediately obtain the slicing inequality
$$
\left(\#K\right)^{\frac{n-1}{n}} \leq  O(n^2) \max\limits_{t, b \in \Z^n, b \not =0}\#(K\cap (t+ b^\perp))
$$
for any origin-symmetric convex body $K$ in $\R^n$. Actually,  Freyer and Henk \cite{FH-2021} removed the condition of $K$ being symmetric
in the statement above (the condition cannot be removed in the discrete analogue of Meyer's inequality). Moreover, they were able to prove
that in the symmetric case the $O(n^2)$-term can be replaced by $O(n)$.  Unfortunately, these results do not seem to apply directly
to Question~1 due to the lack of a direct analogue of Brunn's inequality.

\footnotesize
\bibliographystyle{amsplain}

\end{document}